# Beamforming for Multiuser Massive MIMO Systems: Digital versus Hybrid Analog-Digital

Tadilo Endeshaw Bogale and Long Bao Le
Institute National de la Recherche Scientifique (INRS)
Université du Québec, Montréal, Canada
Email: {tadilo.bogale, long.le}@emt.inrs.ca

*Abstract*— This paper designs a novel hybrid (a mixture of analog and digital) beamforming and examines the relation between the hybrid and digital beamformings for downlink multiuser massive multiple input multiple output (MIMO) systems. We assume that perfect channel state information is available only at the transmitter and we consider the total sum rate maximization problem. For this problem, the hybrid beamforming is designed indirectly by considering a weighed sum mean square error (WSMSE) minimization problem incorporating the solution of digital beamforming which is obtained from the block diagonalization technique. The resulting WSMSE problem is solved by applying the theory of compressed sensing. The relation between the hybrid and digital beamformings is studied numerically by varying different parameters, such as the number of radio frequency (RF) chains, analog to digital converters (ADCs) and multiplexed symbols. Computer simulations reveal that for the given number of RF chains and ADCs, the performance gap between digital and hybrid beamformings can be decreased by decreasing the number of multiplexed symbols. Moreover, for the given number of multiplexed symbols, increasing the number of RF chains and ADCs will increase the total sum rate of the hybrid beamforming which is expected.

## I. INTRODUCTION

Multiple input multiple output (MIMO) system is one of the promising techniques for exploiting the spectral efficiency of wireless channels [1]–[4]. To exploit the full potential of a MIMO system, one must leverage the conventional beamforming techniques (i.e., digital beamforming). Recently the deployment of large antenna arrays at the transmitter and (or) receiver (massive MIMO) has received a lot of attention [5], [6]. Using the law of large numbers, in a rich scattering environment, the works of [5] exploit the fact that the full potential of massive MIMO systems can be achieved by utilizing simple digital beamforming techniques such as zero forcing (ZF) and maximum ratio transmission (MRT).

Despite the research efforts to deploy efficient wireless technologies, wireless industries always face spectrum crunch at microwave frequencies (i.e., frequencies up to 6GHz). Due to this fact, there is an interest to exploit the underutilized millimeter wave (mmWave) frequencies (typical values are $30 - 300$GHz) for cellular applications [7], [8]. However, the mmWave frequencies face severe path loss, penetration loss and rain fading, and they are easily absorbed or scattered by gases [9]. To address these challenges, mmWave frequencies require very high gain antenna systems. It is known that mmWaves frequencies have already used for outdoor point-to-point backhaul links of cellular systems. In the existing backhaul, however, large physical aperture antenna is used to achieve the required link gain which is not economically advantageous. Due to this fact, analog beamforming leveraging massive MIMO is recently suggested for mmWave frequency applications [7]. The fundamental idea of analog beamforming is to control the phase of each antenna's transmitted signal using low cost phase shifters (i.e., each of the analog beamforming coefficients has constant modulus) [7], [10], [11].

At the transmitter side, digital beamfoming is performed at the baseband stage (i.e., the transmitted signals phase and amplitude is determined at baseband frequency). This baseband processing requires dedicated radio frequency (RF) chain (which is expensive) per each antenna. At the receiver side digital beamforming is realized as follows. First, the received signal of each antenna is separately acquired in digital form. Then, the received samples of each antenna are processed jointly to decode the transmitted bits. Thus, in a digital beamforming, the number of received antennas is the same as the number of analog to digital converters (ADCs) which is one of the most expensive parts of a digital receiver [12]. Thus, realizing digital beamforming for massive MIMO systems appears to be impractical. As we have explained previously, the analog beamforming is implemented just by employing phase shifters which are low cost. Therefore, analog beamforming is economically more attractive than that of the digital one. However, as the amplitude of a phase shifter is constant, the performance of analog beamforming is inferior to that of the digital one. To achieve better performance, a hybrid (i.e., a mixture of analog and digital) beamforming is suggested in [8]. In this paper, the hybrid beamforming is designed for single user massive MIMO systems. Furthermore, the work of this paper considers capacity maximization problem to design the hybrid beamforming. However, it is not clear how the hybrid beamforming approach of [8] can be extended to any other design criteria for both single user and multiuser massive MIMO systems.

In the current paper, a novel hybrid beamforming is designed for downlink multiuser massive MIMO systems. It is assumed that perfect channel state information is available only at the transmitter. This could be achieved by simple time division duplex (TDD) training method as in [6]. Obviously the goal of the hybrid beamforming is to achieve very close performance to that of the digital beamforming. Furthermore, it is also reasonable to design the hybrid beamforming in a way that the performance gap between the hybrid and digital beamformings is almost constant for each symbol. This motivates us to design the hybrid beamforming indirectly by considering a weighed sum mean square error (WSMSE) minimization problem while utilizing the solution of digital beamforming. As will be clear later, the mean square error (MSE) weight of each symbol depends on the design criteria and the digital beamforming approach. In this paper, we consider the total sum rate (i.e., design criteria) maximization problem. For this problem, the digital beamforming is designed by applying the well known block diagonalization (BD) approach [13], [14] whereas, the hybrid beamforming is designed by considering a WSMSE minimization problem where the MSE weight of each symbol is selected as the inverse of the square of its digital beamforming gain. The resulting WSMSE problem is solved by leveraging the compressed sensing theory.

After the hybrid bemforming is designed, the relation between hybrid and digital beamformings is studied numerically by varying different parameters, such as the number of multiplexed symbols, RF chains and ADCs. Computer simulations reveal that for the given number of RF chains and ADCs, the performance gap between the digital and hybrid beamformings can be decreased by decreasing the number of multiplexed symbols. Moreover, for the given number of multiplexed symbols, we have noticed that increasing the number of RF chains and ADCs will increase the total sum rate of the hybrid beamforming which is expected.

This paper is organized as follows. Sections II and III discuss the system and channel models. In Section IV, the proposed hybrid beamforming design is presented. In Section V, computer simulations are used to examine the performance of the proposed hybrid beamforming and to study the relation between the hybrid and digital beamformings. Conclusions are drawn in Section VI.

*Notations:* In this paper, upper/lower-case boldface letters denote matrices/column vectors. $\mathbf{X}_{(i,j)}$, $||\mathbf{X}||_F$, $\mathrm{tr}(\mathbf{X})$, $\mathbf{X}^T$, $\mathbf{X}^H$ and $\mathrm{E}(\mathbf{X})$ denote the $(i,j)$th element, Frobenius norm, trace, transpose, conjugate transpose and expected value of $\mathbf{X}$, respectively. $\mathbf{I}_n$ is the identity matrix of size $n \times n$, $\mathbb{C}^{M \times M}$ and $\Re^{M \times M}$ represent spaces of $M \times M$ matrices with complex and real entries, respectively. The acronym null, s.t and i.i.d denote "(right) null space", "subject to" and "independent and identically distributed", respectively.

## II. SYSTEM MODEL

In this section, the digital and hybrid downlink multiuser MIMO system models are discussed. The transmitter equipped with $N$ transmit antennas is serving $K$ users. Each user has $M_k$ antennas to multiplex $S_k$ symbols. The total number of receiver antennas and symbols are $M = \sum_{k=1}^{K} M_k$ and $S = \sum_{k=1}^{K} S_k$, respectively. The entire symbol can be written in a data vector $\mathbf{d} = [\mathbf{d}_1^T, \cdots, \mathbf{d}_K^T]^T$, where $\mathbf{d}_k \in \mathbb{C}^{S_k \times 1}$ is the symbol vector for the $k$th receiver.

### A. Digital Downlink MIMO System Model

In the digital downlink MIMO system, the transmitter precodes $\mathbf{d} \in \mathbb{C}^{S \times 1}$ by using the overall precoder matrix $\mathbf{B} = [\mathbf{B}_1, \cdots, \mathbf{B}_K]$, where $\mathbf{B}_k \in \mathbb{C}^{N \times S_k}$ is the precoder matrix for the $k$th user. The $k$th receiver uses a linear receiver $\mathbf{W}_k \in \mathbb{C}^{M_k \times S_k}$ to recover its symbol $\mathbf{d}_k$ as

$$\hat{\mathbf{d}}_k^D = \mathbf{W}_k^H(\mathbf{H}_k^H \mathbf{B}\mathbf{d} + \mathbf{n}_k) = \mathbf{W}_k^H(\mathbf{H}_k^H \sum_{i=1}^{K} \mathbf{B}_i \mathbf{d}_i + \mathbf{n}_k) \quad (1)$$

where $\mathbf{n}_k \in \mathbb{C}^{M_k \times 1}$ is the additive Gaussian noise at the $k$th receiver and $\mathbf{H}_k^H \in \mathbb{C}^{M_k \times N}$ is the MIMO channel between the transmitter and $k$th receiver and, $\mathbf{B}_k$ and $\mathbf{W}_k, \forall k$ are the conventional digital transmitter (precoder) and receiver (decoder) matrices, respectively. Without loss of generality, we can assume that the entries of $\mathbf{d}_k$ are i.i.d zero mean circularly symmetric complex Gaussian (ZMCSCG) random variables all with unit variance, i.e., $\mathrm{E}\{\mathbf{d}_k \mathbf{d}_k^H\} = \mathbf{I}_{S_k}$, $\mathrm{E}\{\mathbf{d}_k \mathbf{d}_i^H\} = \mathbf{0}$, $\forall i \neq k$, and $\mathrm{E}\{\mathbf{d}_k \mathbf{n}_i^H\} = \mathbf{0}$. The noise vector of the $k$th receiver is a ZMCSCG random variable with variance $\sigma^2$.

### B. Hybrid Downlink MIMO System Model

As discussed in the introduction section of this paper, the hybrid beamforming consists of both digital and analog precoders and decoders. By taking into account of this fact and extending the hybrid beamforming representation of [8] to multiuser case, the estimated signal of the $k$th user can be expressed as

$$\hat{\mathbf{d}}_k^{Hy} = \tilde{\mathbf{W}}_k^H \mathbf{F}_k^H(\mathbf{H}_k^H \sum_{i=1}^{K} \mathbf{A}\tilde{\mathbf{B}}_i \mathbf{d}_i + \mathbf{n}_k) \quad (2)$$

where $\mathbf{A} \in \mathbb{C}^{N \times P_t}$ and $\mathbf{F}_k \in \mathbb{C}^{M_k \times P_{rk}}, \forall k$ are the RF transmitter and receiver matrices (analog beamforming matrices), respectively, and $\tilde{\mathbf{B}}_k$ and $\tilde{\mathbf{W}}_k, \forall k$ are base band (BB) transmitter and receiver matrices (digital beamforming matrices), respectively. Each element of the matrices $\mathbf{A}(\mathbf{F}_k)$ has a constant modulus.

From this equation and (1), we can observe that the output of the digital beamforming will be the same as that of (2) (i.e., $\hat{\mathbf{d}}_k^{Hy} = \hat{\mathbf{d}}_k^D$) when $P_t = N, P_{rk} = M_k, \mathbf{A} = \mathbf{I}_N$ and $\mathbf{F}_k = \mathbf{I}_{M_k}$. Furthermore, (2) becomes the output of analog beamforming when $P_t = S, P_{rk} = S_k, \tilde{\mathbf{B}}_k = \mathbf{I}_{S_k}$ and $\tilde{\mathbf{W}}_k = \mathbf{I}_{S_k}$. From (2), we can also understand that $P_t$ represents the number of transmitter RF chains and $P_{rk}$ represents the number of ADCs at the $k$th receiver. This clearly shows that the performance of $\hat{\mathbf{d}}_k^{Hy}$ depends on the selection of $P_t$ and $P_{rk}$.

## III. CHANNEL MODEL

This section discusses the channel model used in this paper. We consider the most widely used geometric channel model where the channel between the transmitter and the $k$th receiver has $L_k$ scatterers. Under this assumption, $\mathbf{H}_k$ can be expressed as [8], [15]

$$\mathbf{H}_k = \sqrt{\frac{NM_k}{L_k \rho_k}} \sum_{i=1}^{L_k} g_{ki} \mathbf{a}_{tk}(\theta_{tk}(i)) \mathbf{a}_{rk}^H(\theta_{rk}(i)) \quad (3)$$

where $g_{ki}$ is the complex gain of the $k$th user $i$th path with $E\{|g_{ki}|^2\} = 1$, $\rho_k$ is the pathloss between the transmitter and $k$th receiver, $\theta_t(i) \in [0, 2\pi]$, $\theta_{rk}(i) \in [0, 2\pi], \forall i$, and $\mathbf{a}_t(.)$ and $\mathbf{a}_r(.)$ are the antenna array response vectors at the transmitter and receiver, respectively. In particular, this paper adopts a uniform linear array (ULA), where $\mathbf{a}_{tk}(.)$ and $\mathbf{a}_{rk}(.)$ are modeled as [7]

$$\mathbf{a}_{tk}(\theta) = \frac{1}{\sqrt{N}}[1, \exp^{j\frac{2\pi}{\lambda}d\sin(\theta)}, \cdots, \exp^{j(N-1)\frac{2\pi}{\lambda}d\sin(\theta)}]^T$$

$$\mathbf{a}_{rk}(\theta) = \frac{1}{\sqrt{M_k}}[1, \exp^{j\frac{2\pi}{\lambda}d\sin(\theta)}, \cdots, \exp^{j(M_k-1)\frac{2\pi}{\lambda}d\sin(\theta)}]^T$$

where $j = \sqrt{-1}$, $\lambda$ is the transmission wave length and $d$ is the antenna spacing.

## IV. Hybrid Beamforming Design

In this section we design the proposed hybrid beamforming. As we can see from (2), the hybrid beamforming matrices are more constrained than those of the digital ones. Therefore, for any design criteria, the best performance is achieved by digital beamforming. Due to this reason, we utilize $\hat{\mathbf{d}}_k^D$ as a reference received signal. And we quantify the quality of the proposed hybrid beamforming by examining the euclidean distance between the received signal of the hybrid beamforming and that of the digital beamforming. Mathematically, this distance can be obtained by evaluating the MSE between $\hat{\mathbf{d}}_k^{Hy}$ and $\hat{\mathbf{d}}_k^D$ which is given by

$$\begin{aligned}\boldsymbol{\xi}_k =& E\{(\hat{\mathbf{d}}_k^{Hy} - \hat{\mathbf{d}}_k^D)(\hat{\mathbf{d}}_k^{Hy} - \hat{\mathbf{d}}_k^D)^H\} \\ =& E\{(\tilde{\mathbf{W}}_k^H \mathbf{F}_k^H \mathbf{H}_k^H \mathbf{A}\tilde{\mathbf{B}}_k - \mathbf{R}_{kk})\mathbf{d}_k \times \\ & \mathbf{d}_k^H (\tilde{\mathbf{W}}_k^H \mathbf{F}_k^H \mathbf{H}_k^H \mathbf{A}\tilde{\mathbf{B}}_k - \mathbf{R}_{kk})^H + \\ & \sum_{i=1,i\neq k}^K (\tilde{\mathbf{W}}_k^H \mathbf{F}_k^H \mathbf{H}_k^H \mathbf{A}\tilde{\mathbf{B}}_i - \mathbf{R}_{ki})\mathbf{d}_i \times \\ & \mathbf{d}_i^H (\tilde{\mathbf{W}}_k^H \mathbf{F}_k^H \mathbf{H}_k^H \mathbf{A}\tilde{\mathbf{B}}_i - \mathbf{R}_{ki})^H \\ & + (\tilde{\mathbf{W}}_k^H \mathbf{F}_k^H - \mathbf{W}_k^H)\mathbf{n}_k \mathbf{n}_k^H(\tilde{\mathbf{W}}_k^H \mathbf{F}_k^H - \mathbf{W}_k^H)^H\} \\ =& (\tilde{\mathbf{W}}_k^H \mathbf{F}_k^H \mathbf{H}_k^H \mathbf{A}\tilde{\mathbf{B}}_k - \mathbf{R}_{kk}) \times \\ & (\tilde{\mathbf{W}}_k^H \mathbf{F}_k^H \mathbf{H}_k^H \mathbf{A}\tilde{\mathbf{B}}_k - \mathbf{R}_{kk})^H + \\ & \sum_{i=1,i\neq k}^K (\tilde{\mathbf{W}}_k^H \mathbf{F}_k^H \mathbf{H}_k^H \mathbf{A}\tilde{\mathbf{B}}_i - \mathbf{R}_{ki}) \times \\ & (\tilde{\mathbf{W}}_k^H \mathbf{F}_k^H \mathbf{H}_k^H \mathbf{A}\tilde{\mathbf{B}}_i - \mathbf{R}_{ki})^H \\ & + \sigma^2 (\tilde{\mathbf{W}}_k^H \mathbf{F}_k^H - \mathbf{W}_k^H)(\tilde{\mathbf{W}}_k^H \mathbf{F}_k^H - \mathbf{W}_k^H)^H \end{aligned} \quad (4)$$

where $\mathbf{R}_{kk} = \mathbf{W}_k^H \mathbf{H}_k^H \mathbf{B}_k$ and $\mathbf{R}_{ki} = \mathbf{W}_k^H \mathbf{H}_k^H \mathbf{B}_i$. As can be seen from this equation, $\boldsymbol{\xi}_k = 0$ when $\mathbf{A} = \mathbf{I}_N$, $\mathbf{F}_k = \mathbf{I}_{S_k}$, $\tilde{\mathbf{B}}_k = \mathbf{B}_k$, $\tilde{\mathbf{W}}_k = \mathbf{W}_k, \forall k$ which is expected.

As we can see, $\boldsymbol{\xi}_k$ depends on $\mathbf{B}_k, \mathbf{W}_k, \forall k$ which consequently depend on the digital beamforming approach and design criteria. In this paper, we utilize the well known BD digital beamforming approach and employ this approach to maximize the total sum rate of the downlink multiuser massive MIMO system. A brief summary of digital BD beamforming approach is presented in the appendix (see Appendix A for the details).

Using the BD digital beamforming approach, the input output relation of (1) can be rewritten as

$$\begin{aligned}\hat{\mathbf{d}}_k^D =& \mathbf{Z}_k \sqrt{\mathbf{Q}}_k \mathbf{d}_k + \tilde{\mathbf{U}}_{hk}^H \mathbf{n}_k \\ \Rightarrow \hat{d}_{ki}^D =& z_{ki}\sqrt{q_{ki}} d_{ki} + \tilde{\mathbf{u}}_{hki}^H \mathbf{n}_k, \quad \forall k, i\end{aligned} \quad (5)$$

where $\tilde{\mathbf{U}}_{hk} \in \mathbb{C}^{M_k \times S_k}$ is a unitary matrix, $\mathbf{Z}_k \in \mathbb{C}^{S_k \times S_k}$ is a diagonal matrix, $\mathbf{Q}_k \in \Re^{S_k \times S_k}$ is a real diagonal non-negative power allocation matrix, $\hat{d}_{ki}^D(d_{ki})$ is the $i$th element of $\hat{\mathbf{d}}_k^D(\mathbf{d}_k)$ and $z_{ki}(q_{ki})$ is the $i$th diagonal element of $\mathbf{Z}_k(\mathbf{Q}_k)$ and $\tilde{\mathbf{u}}_{hki}$ is the $i$th row of $\tilde{\mathbf{U}}_{hk}$.

From this equation, we can observe that in the BD beamforming, $\mathbf{R}_{kk} = \mathbf{Z}_k\sqrt{\mathbf{Q}}_k$ and $\mathbf{R}_{ki} = \mathbf{0}, \forall k \neq i$. Consequently, $\boldsymbol{\xi}_k$ of (4) can be expressed as

$$\begin{aligned}\boldsymbol{\xi}_k =& \\ & (\tilde{\mathbf{W}}_k^H \mathbf{F}_k^H \mathbf{H}_k^H \mathbf{A}\tilde{\mathbf{B}}_k - \mathbf{Z}_k\sqrt{\mathbf{Q}}_k)(\tilde{\mathbf{W}}_k^H \mathbf{F}_k^H \mathbf{H}_k^H \mathbf{A}\tilde{\mathbf{B}}_k - \mathbf{Z}_k\sqrt{\mathbf{Q}}_k)^H \\ & + \sum_{i=1,i\neq k}^K (\tilde{\mathbf{W}}_k^H \mathbf{F}_k^H \mathbf{H}_k^H \mathbf{A}\tilde{\mathbf{B}}_i)(\tilde{\mathbf{W}}_k^H \mathbf{F}_k^H \mathbf{H}_k^H \mathbf{A}\tilde{\mathbf{B}}_i)^H \\ & + \sigma^2(\tilde{\mathbf{W}}_k^H \mathbf{F}_k^H - \tilde{\mathbf{U}}_{hk}^H)(\tilde{\mathbf{W}}_k^H \mathbf{F}_k^H - \tilde{\mathbf{U}}_{hk}^H)^H.\end{aligned} \quad (6)$$

As discussed previously, it is also reasonable to design $\mathbf{A}, \tilde{\mathbf{B}}_k, \tilde{\mathbf{W}}_k$ and $\mathbf{F}_k$ in a way that the performance gap between the hybrid and digital beamformings is almost constant for all symbols. Due to this reason, we suggest a WSMSE optimization problem as our objective function, where the weight of the $k$th user $i$th symbol is set to $\frac{1}{z_{ki}^2 q_{ki}}$. This problem can be mathematically formulated as[1]

$$\min_{\mathbf{A}, \tilde{\mathbf{B}}_k, \tilde{\mathbf{W}}_k, \mathbf{F}_k} \sum_{k=1}^K \text{tr}\{(\mathbf{Z}_k\sqrt{\mathbf{Q}}_k)^{-1}\boldsymbol{\xi}_k(\mathbf{Z}_k\sqrt{\mathbf{Q}}_k)^{-1}\} \triangleq \xi_w$$

$$\text{s.t} \sum_{k=1}^K \text{tr}\{\mathbf{A}\tilde{\mathbf{B}}_k\tilde{\mathbf{B}}_k^H\mathbf{A}^H\} = P_{max}$$

$$|\mathbf{A}_{(i,j)}|^2 = \kappa, |\mathbf{F}_{k(i,j)}|^2 = \kappa_k \quad (7)$$

where the constraint $\sum_{k=1}^K \text{tr}\{\mathbf{A}\tilde{\mathbf{B}}_k\tilde{\mathbf{B}}_k^H\mathbf{A}^H\} = P_{max}$ is introduced to ensure that both the digital and hybrid beamformings utilize the same total power, and $\kappa$ and $\kappa_k, \forall k$ are arbitrary positive constants. After some mathematical manipulation, one can express $\xi_w$ as

$$\xi_w = ||\tilde{\tilde{\mathbf{W}}}^H \mathbf{H}^H \mathbf{A}\tilde{\mathbf{B}} - \mathbf{I}||_F^2 + \sum_{k=1}^K ||\tilde{\tilde{\mathbf{W}}}_k^H \mathbf{F}_k^H - \tilde{\tilde{\mathbf{U}}}_{hk}^H||_F^2 \quad (8)$$

where $\tilde{\mathbf{B}} = [\tilde{\mathbf{B}}_1, \cdots, \tilde{\mathbf{B}}_K]$, $\tilde{\tilde{\mathbf{W}}}_k = \tilde{\mathbf{W}}_k(\mathbf{Z}_k\sqrt{\mathbf{Q}}_k)^{-1}$, $\tilde{\tilde{\mathbf{U}}}_{hk} = \tilde{\mathbf{U}}_{hk}(\mathbf{Z}_k\sqrt{\mathbf{Q}}_k)^{-1}$, $\tilde{\tilde{\mathbf{W}}} = \text{blkdiag}(\tilde{\tilde{\mathbf{W}}}_1, \cdots, \tilde{\tilde{\mathbf{W}}}_K)$ and $\mathbf{F} = \text{blkdiag}(\mathbf{F}_1, \cdots, \mathbf{F}_K)$. The above optimization problem can thus be re-expressed as

$$\min_{\mathbf{A}, \tilde{\mathbf{B}}, \tilde{\tilde{\mathbf{W}}}_k, \mathbf{F}_k} ||(\mathbf{F}\tilde{\tilde{\mathbf{W}}})^H \mathbf{H}^H \mathbf{A}\tilde{\mathbf{B}} - \mathbf{I}||_F^2 + \sum_{k=1}^K ||\tilde{\tilde{\mathbf{W}}}_k^H \mathbf{F}_k^H - \tilde{\tilde{\mathbf{U}}}_{hk}^H||_F^2$$

$$\text{s.t} \quad \text{tr}\{\mathbf{A}\tilde{\mathbf{B}}\tilde{\mathbf{B}}^H \mathbf{A}^H\} = P_{max}$$

$$|\mathbf{A}_{(i,j)}|^2 = \kappa, |\mathbf{F}_{k(i,j)}|^2 = \kappa_k. \quad (9)$$

---
[1]Note that when $z_{ki}^2 q_{ki} \approx 0$, we ignore the $k$th user $i$th symbol. This is because, in such a case, the digital beamforming also switches off this symbol.

To solve this problem, we employ two steps. In the first step, we jointly optimize $\tilde{\tilde{\mathbf{W}}}_k$ and $\mathbf{F}_k$ for each user. In the second step, we optimize $\mathbf{A}$ and $\tilde{\mathbf{B}}$ jointly for fixed $\tilde{\tilde{\mathbf{W}}}_k$ and $\mathbf{F}_k, \forall k$. In the following, we provide detailed explanation of these two steps

## A. Step 1

For the $k$th user, the joint optimization of $\tilde{\tilde{\mathbf{W}}}_k, \mathbf{F}_k$ can be expressed as

$$\min_{\tilde{\tilde{\mathbf{W}}}_k, \mathbf{F}_k} ||\mathbf{F}_k \tilde{\tilde{\mathbf{W}}}_k - \tilde{\mathbf{U}}_{hk}||_F^2, \quad \text{s.t } |\mathbf{F}_{k(i,j)}|^2 = 1 \quad (10)$$

where $\kappa_k$ is set to unity without loss of generality. As we can see, $\tilde{\tilde{\mathbf{W}}}_k$ and $\mathbf{F}_k$ are jointly coupled in the objective function. Furthermore, the constraint function is non convex. So the global optimal solution of this problem is not anticipated. Here our aim is to provide suboptimal (close to optimal) solution. To this end, we leverage the technique of compressed sensing approach to solve this problem which is presented as follows.

As we can see from (17) (i.e., in Appendix A), $\tilde{\mathbf{U}}_{hk}$ is highly correlated with the left singular eigenvectors of the channel $\mathbf{H}_k^H$. These eigenvectors are likely formed from the linear combinations of $\tilde{\mathbf{F}}_k$ and the null spaces of $\tilde{\mathbf{F}}_k^T$ where $\tilde{\mathbf{F}}_k = [\mathbf{a}_{rk}(\theta_{rk}(1)), \mathbf{a}_{rk}(\theta_{rk}(2)), \cdots, \mathbf{a}_{rk}(\theta_{rk}(L_k))]$. Due to this intuition, we assume that each column of $\mathbf{F}_k$ is taken from one of the columns of $\mathbb{F}_k$, where $\mathbb{F}_k = [\tilde{\mathbf{F}}_k \text{ null}(\tilde{\mathbf{F}}_k^T)] \in \mathbb{C}^{M_k \times M_k}$ with normalized entries (i.e., the modulus of each element of $\mathbb{F}_k$ is normalized to 1). By applying this idea, we reexpress the above problem as

$$\min_{\tilde{\tilde{\mathbf{W}}}_k, \mathbf{F}_k} ||\mathbf{F}_k \tilde{\tilde{\mathbf{W}}}_k - \tilde{\mathbf{U}}_{hk}||_F^2 \quad \text{s.t } \mathbf{F}_{k(:,i)} \in \mathbb{F}_{k(:,\forall j)}. \quad (11)$$

This problem is still non convex. To simplify this problem, we introduce the following variable $\bar{\mathbf{F}}_\mathbf{k} = \mathbb{F}_k$. Using this variable, we can simplify this problem as

$$\min_{\bar{\mathbf{W}}_k} ||\bar{\mathbf{F}}_k \bar{\mathbf{W}}_k - \tilde{\mathbf{U}}_{hk}||_F^2 \quad \text{s.t } ||\text{diag}(\bar{\mathbf{W}}_k \bar{\mathbf{W}}_k^H)||_0 = P_{rk} \quad (12)$$

where the equality constraint is introduced to ensure that the dimension of $\bar{\mathbf{W}}_k$ is the same as that of $\tilde{\tilde{\mathbf{W}}}_k$. This problem is a dimension reduction problem. As this problem contains non-convex constraint, convex optimization can not be applied. However, plenty of compressed sensing algorithms can be applied to solve this problem. In the current paper, we utilize simple orthogonal matching pursuit approach to solve the problem [16]–[18]. Orthogonal matching pursuit is a numerical approach of finding the best matching projections of multidimensional data onto an over-complete dictionary (for our problem $\bar{\mathbf{F}}_k$ is the dictionary matrix of the $k$th user). The detailed explanation of the matching pursuit algorithm can be found in [16].

For our problem, the orthogonal matching pursuit algorithm is summarized as follows.
**Algorithm I**: Orthogonal matching pursuit to solve (12)

1) **Initialization**: Set $\mathbf{F}_k = [\,]$, $\bar{\mathbf{F}}_{kR} = \tilde{\mathbf{U}}_{hk}$.
2) **for** i=1:$P_{rk}$ **do**
   $\boldsymbol{\Phi} = \bar{\mathbf{F}}_k^H \bar{\mathbf{F}}_{kR}$.
   $i = \arg\max(\text{diag}(\boldsymbol{\Phi}\boldsymbol{\Phi}^H))$.
   $\mathbf{F}_k = [\mathbf{F}_k | \bar{\mathbf{F}}_{(:,i)}]$.
   $\tilde{\tilde{\mathbf{W}}}_k = (\mathbf{F}_k^H \mathbf{F}_k)^{-1} \mathbf{F}_k^H \tilde{\mathbf{U}}_{hk}$.
   $\bar{\mathbf{F}}_{kR} = \frac{\tilde{\mathbf{U}}_{hk} - \mathbf{F}_k \tilde{\tilde{\mathbf{W}}}_k}{||\tilde{\mathbf{U}}_{hk} - \mathbf{F}_k \tilde{\tilde{\mathbf{W}}}_k||_F}$.
   **end for**

## B. Step 2

For the given $\tilde{\tilde{\mathbf{W}}}_k$ and $\tilde{\mathbf{F}}_k, \forall k$, (9) can be simplified to

$$\min_{\mathbf{A}, \tilde{\mathbf{B}}} ||\tilde{\tilde{\mathbf{W}}}^H \mathbf{F}^H \mathbf{H}^H \mathbf{A} \tilde{\mathbf{B}} - \mathbf{I}||_F^2$$
$$\text{s.t } \text{tr}\{\mathbf{A}\tilde{\mathbf{B}}\tilde{\mathbf{B}}^H \mathbf{A}^H\} \leq P_{max}$$
$$|\mathbf{A}_{(i,j)}|^2 = 1. \quad (13)$$

From (17) (see Appendix A), we can also notice that $\mathbf{A}$ is highly correlated with the right singular eigenvectors of $\mathbf{H}^H$. These eigenvectors are likely formed from the linear combinations of $\tilde{\mathbf{A}}$ and the null spaces of $\tilde{\mathbf{A}}^T$, where $\tilde{\mathbf{A}} = [\mathbf{a}_{tk}(\theta_{t1}(1)), \cdots, \mathbf{a}_{tk}(\theta_{t1}(L_1)), \cdots, \mathbf{a}_{tk}(\theta_{tK}(1)), \cdots, \mathbf{a}_{tk}(\theta_{tK}(L_K))]$. Because of this intuition, it is assumed that each column of $\mathbf{A}$ is taken from one of the columns of $\mathbb{A}$, where $\mathbb{A} = [\tilde{\mathbf{A}} \text{ null}(\tilde{\mathbf{A}}^T)] \in \mathbb{C}^{N \times N}$ with normalized entries (i.e., the modulus of each element of $\mathbb{A}$ is normalized to 1). By employing $\bar{\mathbf{A}} = \mathbb{A}$ and compressed sensing technique in Step 1, the above problem can be reexpressed as

$$\min_{\bar{\mathbf{B}}} ||(\mathbf{H}\mathbf{F}\tilde{\tilde{\mathbf{W}}})^H \bar{\mathbf{A}}\bar{\mathbf{B}} - \mathbf{I}||_F^2$$
$$\text{s.t } \text{tr}\{\bar{\mathbf{A}}\bar{\mathbf{B}}\bar{\mathbf{B}}^H \bar{\mathbf{A}}^H\} = P_{max}$$
$$||\text{diag}(\bar{\mathbf{B}}\bar{\mathbf{B}}^H)||_0 = P_t. \quad (14)$$

This problem can be solved by employing the matching pursuit algorithm like that of (12) and is summarized as follows.
**Algorithm II**: Orthogonal matching pursuit to solve (14)

1) **Initialization**: Set $\mathbf{A} = [\,]$, $\bar{\mathbf{A}}_R = \mathbf{I}$, $\bar{\bar{\mathbf{A}}} = (\mathbf{H}\mathbf{F}\tilde{\tilde{\mathbf{W}}})^H \bar{\mathbf{A}}$.
2) **for** i=1:$P_t$ **do**
   $\boldsymbol{\Phi} = \bar{\bar{\mathbf{A}}}^H \bar{\mathbf{A}}_R$.
   $i = \arg\max(\text{diag}(\boldsymbol{\Phi}\boldsymbol{\Phi}^H))$.
   $\mathbf{A} = [\mathbf{A}|\bar{\mathbf{A}}_{(:,i)}]$.
   $\hat{\mathbf{A}} = (\mathbf{H}\mathbf{F}\tilde{\tilde{\mathbf{W}}})^H \mathbf{A}$.
   $\tilde{\mathbf{B}} = (\hat{\mathbf{A}}^H \hat{\mathbf{A}})^{-1} \hat{\mathbf{A}}^H \mathbf{I}$.
   $\bar{\mathbf{A}}_R = \frac{\mathbf{I} - \hat{\mathbf{A}}\tilde{\mathbf{B}}}{||\mathbf{I} - \hat{\mathbf{A}}\tilde{\mathbf{B}}||_F}$.
   **end for**
3) $\tilde{\tilde{\mathbf{B}}} = \sqrt{P_{max}} \frac{\tilde{\mathbf{B}}}{||\mathbf{A}\tilde{\mathbf{B}}||_F}$.

We would like to mention here that in practice the optimization problems (12) and (14) are solved at the transmitter. However, the $k$th receiver require $\mathbf{F}_k$ (i.e., $M_k \times P_{rk}$) and $\tilde{\tilde{\mathbf{W}}}_k$ (i.e., $P_{rk} \times S_k$) to decode its own data. As each column of $\mathbf{F}_k$ is selected from the columns of $\mathbb{F}_k$, $\mathbf{F}_k$ can be constructed from $\theta_{rk}(i)$. From this explanation, we can notice that the $k$th receiver can locally compute $\mathbf{F}_k$ just by receiving $\theta_{rk}(i)$ (i.e., $L_k$ real scalars) and the column index of $\mathbb{F}_k$ corresponding to $\mathbf{F}_k$ (i.e., $P_{rk}$ real scalars). This shows that the feedback information from the transmitter to the $k$th receiver is $L_k + P_{rk}$ real plus $P_{rk} \times S_k$ complex scalars which is insignificant for practically relevant scenario (Typically $L_k$ is in the order of 10 to 20, and $P_k$ and $S_k$ are in the orders of 1 to 10).

## V. SIMULATION RESULTS

This section presents simulation results. We have used $N = 128$, $M_k = 32$, $d = 0.5\lambda$, $K = 4$, $L_k = 16$ and $P_t = KP_{rk}$. The signal to noise ratio (SNR) which is defined as $SNR = \frac{P_{av}}{\sigma^2}$ is controlled by varying $\sigma^2$ while keeping the total transmitted power $P_{max} = KS_k$ mw and $P_{av} = \frac{P_{max}}{S_k}$. The channel parameters $\boldsymbol{\rho} = [0.2338\ 0.2333\ 0.0402\ 0.5290]$ and $\theta_{tk}(i)(\theta_{rk}(i))$ are taken randomly from uniform random variables in $[0,\ 2\pi]$ (see (15)). The total sum rate is computed as $R_t = \sum_{k=1}^{K}\sum_{i=1}^{S_k} R_{ki}$, where the rate of each symbol is calculated as $R_{ki} = \log_2(1 + \gamma_{ki})$ with $\gamma_{ki}$ as the achieved signal to interference plus noise ratio (SINR) of the $k$th user $i$th symbol. All of the plots are generated by averaging 1000 realizations of $g_k, \forall k$.

$$\boldsymbol{\theta}_t(\text{in } \pi) = \begin{bmatrix} 0.71 & 0.66 & 0.20 & 0.70 \\ 1.29 & 0.16 & 0.81 & 1.88 \\ 0.65 & 0.91 & 0.22 & 0.52 \\ 1.70 & 0.64 & 0.94 & 1.61 \\ 0.80 & 1.85 & 1.42 & 1.33 \\ 0.43 & 1.11 & 1.99 & 1.38 \\ 1.24 & 0.61 & 1.53 & 1.95 \\ 0.80 & 0.10 & 0.77 & 1.22 \\ 0.47 & 0.07 & 0.38 & 0.39 \\ 0.43 & 1.50 & 0.83 & 1.83 \\ 0.04 & 0.71 & 1.27 & 1.24 \\ 0.60 & 0.53 & 0.16 & 0.32 \\ 1.50 & 1.76 & 0.08 & 0.71 \\ 0.77 & 1.48 & 0.05 & 1.23 \\ 0.74 & 1.89 & 1.70 & 0.20 \\ 0.02 & 0.83 & 0.36 & 1.71 \end{bmatrix} \quad \boldsymbol{\theta}_r(\text{in } \pi) = \begin{bmatrix} 0.25 & 1.87 & 0.73 & 0.45 \\ 1.04 & 0.73 & 0.73 & 0.16 \\ 0.39 & 0.74 & 1.51 & 1.30 \\ 0.23 & 1.52 & 1.69 & 1.40 \\ 0.77 & 0.16 & 0.55 & 1.04 \\ 1.92 & 1.63 & 1.81 & 1.71 \\ 1.80 & 1.58 & 0.07 & 1.86 \\ 1.37 & 0.44 & 0.93 & 0.16 \\ 1.85 & 1.05 & 1.19 & 0.75 \\ 0.49 & 0.64 & 1.31 & 0.35 \\ 0.20 & 1.09 & 1.28 & 0.24 \\ 0.85 & 1.39 & 0.98 & 1.20 \\ 1.13 & 0.10 & 1.13 & 1.46 \\ 0.88 & 1.92 & 1.59 & 1.08 \\ 0.84 & 1.93 & 1.68 & 1.47 \\ 1.56 & 0.25 & 1.44 & 1.96 \end{bmatrix} \quad (15)$$

### A. Comparison of Digital and Hybrid Beamformings

In this subsection, we compare the performance of the digital and hybrid beamforming algorithms. To this end, we set $S_k = 8$ and $P_{rk} = 16, \forall k$. Fig. 1 shows the sum rate achieved by the digital and hybrid beamformings. As can be seen from this figure, the performance of the hybrid beamforming is almost the same as that of the digital one from the low to moderate SNR regions. And small performance gap is observed at high SNR regions.

### B. Effect of $P_{rk}$ on Hybrid Beamforming

In this subsection, we examine the effect of $P_{rk}$ on the performance of the hybrid beamforming algorithm. Fig. 2 illustrates the performance of the hybrid beamforming for different settings of $P_{rk}$ at SNR=$\{-4\text{dB}, 6\text{dB}\}$ and $S_k = 8$. From this figure, one can observe that the performance of the hybrid beamforming algorithm improves as the number of RF chains and ADCs (i.e., $P_{rk}$) increase which is expected.

### C. Joint effects of $S_k$ and $P_{rk}$ on Digital and Hybrid Beamformings

In this simulation, we examine the joint effects of $S_k$ and $P_{rk}$ on the performance of the digital and hybrid beamformings. To this end, we employ SNR=$-4$dB. Fig. 3 shows the performance of the digital and hybrid beamformings. As can be seen from this figure, the achieved sum rates of the digital and hybrid beamformings increase as $S_k$ increases. This is because as $S_k$ increases, $P_{max} = KS_k$ also increases. Also the performance gap between hybrid and digital beamformings decrease as $S_K$ decreases. Hence, for a limited number of RF chains and ADCs (i.e., $P_{rk}$), the performance gap between hybrid and digital beamformings can be decreased by reducing the number of multiplexed symbols ($S_k$).

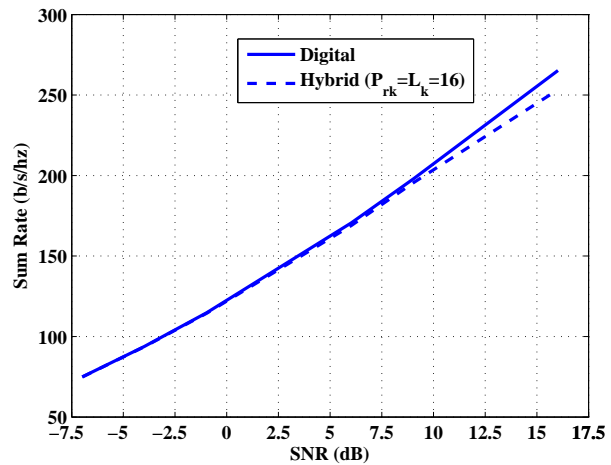

Fig. 1. Comparison of digital and hybrid beamformings when $L_k = 16$, $S_k = 8$ and $P_{rk} = 16$.

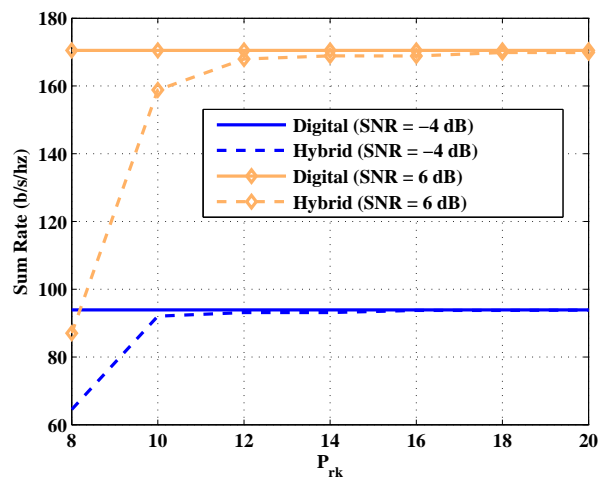

Fig. 2. Effect of $P_{rk}$ on the hybrid beamforming.

## VI. CONCLUSIONS

This paper designs novel hybrid beamforming and discusses the relation between hybrid and digital beamformings for downlink multiuser massive MIMO systems. The design and analysis is provided by considering the total sum rate maximization problem. For this problem, the hybrid beamforming is designed indirectly by considering a WSMSE minimization problem incorporating the solution of digital beamforming which is obtained from the BD technique. The resulting WSMSE problem is solved by leveraging the compressed sensing theory. After the hybrid bemforming is designed, the relation between hybrid and digital beamformings is studied numerically by varying different parameters, such as the number of multiplexed symbols, RF chains and ADCs. Simulation results show that, for the given number of RF chains and ADCs, the performance gap between the digital and hybrid beamformings can be decreased by decreasing the number of multiplexed symbols. Moreover, for the given number of multiplexed symbols, increasing the number of RF chains and ADCs increases the achievable sum rate of the hybrid beamforming which is expected.

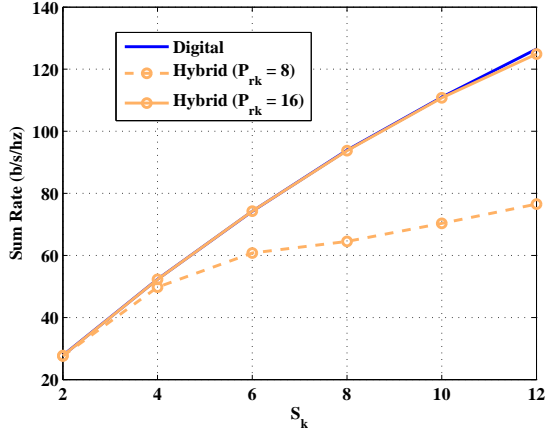

Fig. 3. Effect of $S_k$ and $P_{rk}$ on the digital and hybrid beamforming.

## APPENDIX A
## BLOCK DIAGONALIZATION BEAMFORMING

In this appendix, the well known digital BD beamforming algorithm for multiuser MIMO systems is summarized. For each desired user, the BD beamforming algorithm utilizes three key steps. In the first step, the interference of the other users is eliminated either partially or completely (if possible). In the second step, the self interference of each user is canceled (i.e., each of the symbols of the desired user are parallel). In the third step, the powers of each symbol is optimized to maximize the total sum rate of the downlink system [13], [14].

For better exposition, let us define the following terms
$$\tilde{\mathbf{H}}_k \triangleq [\mathbf{H}_1, \mathbf{H}_2, \cdots, \mathbf{H}_{k-1}\mathbf{H}_{k+1}, \cdots, \mathbf{H}_K]$$
$$\mathbf{V}_{h0k} \triangleq \text{null}(\tilde{\mathbf{H}}_k^H) \quad (16)$$
$$\mathbf{X}_k \triangleq \mathbf{H}_k^H \mathbf{V}_{0hk} = \tilde{\mathbf{U}}_{hk} \mathbf{Z}_k \tilde{\mathbf{V}}_{hk}^H \quad (17)$$

where $\mathbf{V}_{0hk} \in \mathbb{C}^{N \times S_k}$ is the first $S_k$ vectors of $\text{null}(\tilde{\mathbf{H}}_k^H)$, $\tilde{\mathbf{U}}_{hk} \in \mathbb{C}^{M_k \times S_k}$, $\tilde{\mathbf{V}}_{hk} \in \mathbb{C}^{S_k \times S_k}$ and $\mathbf{Z}_k$ is a diagonal matrix of size $S_k$, and the last equality is derived by applying the singular value decomposition (SVD) operation on $\mathbf{X}_k$.

The first and second steps of the BD beamforming algorithm can be performed by setting $\mathbf{B}_k$ and $\mathbf{W}_k$ of (1) as $\mathbf{B}_k = \mathbf{V}_{0hk}\tilde{\mathbf{V}}_{hk}\sqrt{\mathbf{Q}_k}$ and $\mathbf{W}_k = \tilde{\mathbf{U}}_{hk}$, where $\mathbf{Q}_k$ is the power allocation matrix of the $k$th user. By doing so, the input output relation of (1) can be rewritten as

$$\hat{\mathbf{d}}_k^D = \mathbf{Z}_k\sqrt{\mathbf{Q}_k}\mathbf{d}_k + \tilde{\mathbf{U}}_{hk}^H\mathbf{n}_k, \Rightarrow \hat{d}_{ki}^D = z_{ki}\sqrt{q_{ki}}d_{ki} + \tilde{\mathbf{u}}_{hki}^H\mathbf{n}_k \quad (18)$$

where $\hat{d}_{ki}^D(d_{ki})$ is the $i$th element of $\hat{\mathbf{d}}_k^D(\mathbf{d}_k)$ and $z_{ki}(q_{ki})$ is the $i$th diagonal element of $\mathbf{Z}_k(\mathbf{Q}_k)$ and $\tilde{\mathbf{u}}_{hki}$ is the $i$th row of $\tilde{\mathbf{U}}_{hk}$.

The third step of the BD beamforming is achieved by solving the following sum rate maximization problem

$$\max_{q_{ki}} \sum_{k=1}^{K}\sum_{i=1}^{S_k} \log_2\left(1 + \frac{z_{ki}^2 q_{ki}}{\sigma^2}\right), \text{ s.t } \sum_{k=1}^{K}\sum_{i=1}^{S_k} q_{ki} = P_{max}$$

where $q_{ki}$ is the $i$th diagonal element of $\mathbf{Q}_k$ and $P_{max}$ is the maximum power available at the transmitter. The global optimal solution of this problem can be obtained by simple water filling algorithm [2].